\documentclass[11pt,twoside,reqno]{amsart}
\usepackage{amsmath,amsfonts,amsthm,amssymb}
\usepackage{graphicx,psfrag,overpic}
\usepackage{bm}
\usepackage[all,cmtip,line]{xy}
\usepackage{array,setspace}
\usepackage{color}
\numberwithin{equation}{section}

\usepackage{graphics}
\usepackage{epsfig}
\newcommand{\R}{\mathbb{R}}

\begin{document}

\title[Entropy and Convexity for Nonlinear PDEs]
{Introduction\\[5mm]
Entropy and Convexity\\
 for Nonlinear Partial Differential Equations}

\author{John M. Ball}
\address{John M. Ball, Mathematical Institute, University of Oxford, Oxford, OX2 6GG, UK}
\email{ball@maths.ox.ac.uk}

\author{Gui-Qiang G. Chen}
\address{Gui-Qiang G. Chen, Mathematical Institute, University of Oxford, Oxford, OX2 6GG, UK}
\email{chengq@maths.ox.ac.uk}


\bigskip
\bigskip

\begin{abstract}\noindent Partial differential equations are ubiquitous in almost all applications of mathematics,
where they provide a natural mathematical description of many phenomena involving change in physical,
chemical, biological, and social processes.
The concept of entropy originated in thermodynamics and statistical physics during the $19^{\rm th}$ century
to describe the heat exchanges that occur in the thermal processes in a thermodynamic system,
while the original notion of convexity is for sets and functions in mathematics.
Since then, entropy and convexity have become two of the most important concepts in mathematics.
In particular, nonlinear methods via entropy and convexity have been playing an increasingly important
role in the analysis of nonlinear partial differential equations in recent decades.
This opening article of the theme issue is intended to provide an introduction to
entropy, convexity, and related nonlinear methods for the analysis of nonlinear partial differential
equations.
We also provide a brief discussion about the content and contributions of the papers
that make up this theme issue.
\end{abstract}

\keywords{entropy, convexity, nonlinear methods, partial differential equations,
calculus of variations, hyperbolic conservation laws,
discontinuous solutions, singular solutions,
shock waves, entropy methods, convex analysis,
continuum mechanics, statistical physics, gradient flow, convection, kinetic theory, dynamical systems,
stochastic processes, multi-dimensional dynamics}

\maketitle

\section*{}\label{sec1}

Partial differential equations (PDEs) are ubiquitous in almost all applications of mathematics,
where they provide a natural mathematical description of many phenomena involving change in physical,
chemical, biological, and social processes.
The behaviour of every material object, with length scales ranging from sub-atomic to
astronomical and timescales ranging from picoseconds to millennia,
can be modelled by PDEs or by equations having similar features.

The concept of entropy originated in thermodynamics and statistical physics
during the $19^{\rm th}$ century to describe the heat exchanges that occur in the thermal
processes in a thermodynamic system,
while the original notion of convexity is for sets and functions
in mathematics.
Since then, entropy and convexity have become two of the most important
concepts in mathematics.
In particular, nonlinear methods via entropy and convexity have been playing an
increasingly important role in the analysis of nonlinear PDEs in recent
decades.
For example, for discontinuous or singular solutions to nonlinear conservation
laws which may contain shock waves
and concentrations,
the notion of entropy
solutions is based on entropy conditions involving convexity,
motivated by and consistent with the Second Law of Thermodynamics.
In addition, entropy methods have become one of the most efficient methods
in the analysis
of physically relevant discontinuous or singular solutions.
The notions of convexity appropriate for multi-dimensional problems,
such as polyconvexity, quasiconvexity, and rank-one convexity, are responsible
for several recent major advances.
In the last three decades, various nonlinear methods involving entropy and convexity
have been developed to deal with discontinuous and singular solutions
in different areas of PDEs,
especially in nonlinear conservation laws, the calculus of variations, and
gradient flows.

This theme issue is devoted to fundamental questions concerning
entropy, convexity, and related nonlinear methods designed to
help understand the very difficult problems posed by multi-dimensional,
nonlinear PDE problems.
In particular, it includes the discussion of several recent developments
in nonlinear methods
via entropy and convexity, the exploration of
their underlying connections, and the development of
new unifying methods, ideas, and insights
involving entropy and convexity for important multi-dimensional PDE problems
in fluid/solid mechanics and other areas.
These new developments are at the forefront of current research.

\section{Entropy}
The concept and name of {\it entropy}, as a mathematical quantity,
originated in the early 1850s in the work of
Rudolf Julius Emmanuel Clausius (1822--1888) (see \cite{Clausius}),
built on the previous intuition of
Nicolas L\'{e}onard Sadi Carnot (1796--1832) (see \cite{Carnot});
{\it Entropy}, as an extensive thermodynamic function of state, describes
the heat exchanges that
occur in thermal processes from the macroscopic point of view.
Ludwig Edward Boltzmann (1844--1906) first observed that the Clausius entropy associated with a system in equilibrium
is proportional to the logarithm of the number of microstates in microscopic dynamics
which form the macrostate of this equilibrium; this exposed a symbiotical
relation between the notions of {\it entropy} from both macroscopic and microscopic points
of view (see \cite{Boltzmann}).
Similar ideas were also developed by many physicists and mathematicians of those times,
most notably James Clerk Maxwell (1831--1879), Josiah Willard Gibbs (1839--1903),
Max Karl Ernst Ludwig Planck (1858--1947), as well as Constantin Carath\'{e}odory (1873--1950);
see \cite{Caratheodory,Gibbs,Maxwell,Planck}.
There remain many  issues concerning the relation between the microscopic and macroscopic
descriptions of matter which are not completely resolved,
notably concerning the manner in which irreversibility at the macroscopic level
is consistent with reversible dynamics at the microscopic level.

A key contribution of Boltzmann ({\it cf.} \cite{Boltzmann}) was to notice that the Maxwell-Boltzmann equation
of the kinetic theory of gases possesses a Lyapunov function related to entropy,
the Boltzmann $H$-function.
Macroscopic statements of the Second Law of Thermodynamics, such as the Clausius-Duhem inequality,
deliver similar Lyapunov functions for solutions satisfying the balance laws of mass,
momentum, and energy (see Duhem \cite{Duhem} and Ericksen \cite{Ericksen}),
which, depending on the boundary conditions, is either the total entropy
or the so-called ballistic free energy (when the boundary is in contact with a heat bath
maintained at a temperature
that is constant in space and time). The situation is more complicated for spatially-dependent
boundary temperatures; see [3].
The existence of such Lyapunov functions links dynamics to energy minimization,
providing a dynamic motivation for the calculus of variations.
However, there remains a large theoretical gap between Boltzmann's $H$-theorem,
applying as it does to a moderately rarified gas,
and macroscopic statements of the Second Law,
which are routinely applied to a much wider class of materials.

Since the work of these pioneers, entropy has become one of the most important concepts in the sciences,
and various entropy approaches have been developed which have had a considerable impact
on many important areas
of mathematics and related sciences. These areas include
continuum mechanics, kinetic theory, statistical physics,
probability,
stochastic processes and random fields, dynamical systems and ergodic theory,
information and coding, data analysis and statistical inference,
in addition to PDEs.
Ideas involving entropy have played a crucial role in developing
many important mathematical approaches and methods
such as variational principles, Lyapunov functionals, relative entropy methods,
monotonicity methods, weak convergence methods, divergence-measure fields,
and kinetic methods
({\it cf.}
\cite{BallKnowles,Chen,CF,ChenPerepelitsa,Dacorogna,Dafermos-book,Evans1,Evans2,GSR,Lax}
and \cite{ADPZ}--\cite{Slemrod}).

Entropy methods for nonlinear PDEs
are techniques and approaches
for discovering and exploiting dissipation inequalities for nonlinear PDEs.
As noted in Evans \cite{Evans3}, these approaches reflect the insight
that nonlinear PDEs are generally
too hard to grapple with directly, and so often a good idea is to simplify by integrating out
various expressions involving the solutions to gain more information about them.
The advantage of entropy methods is the elegance in their formulation
and generality in their
applications ({\it cf.} \cite{CF,ChenPerepelitsa,Evans1,Brenier,Dafermos1,Liu,Saint-Raymond}).

\section{Convexity}

The original notion of convexity in mathematics is for sets and functions.
In Euclidean space, a set is convex if, for every pair of points within the set, every point on
the straight line segment that joins them is also within the set.
A real-valued function defined on an interval is said to be convex if the graph of the function
does not lie above the line segment joining any two points of the graph.
It has been generalised to more significant settings, especially for functionals,
which have played an important
role in many areas of mathematics, including functional analysis, complex analysis,
differential geometry, topology, geometric measure theory,
optimization theory, the calculus of variations,
in addition to PDEs.
Convex analysis has become an important branch of mathematics devoted to the study of properties of
convex functions, convex sets, and convex functionals, which occur in the analysis
of PDEs.

In particular, convexity plays an important role for entropy methods for PDEs,
especially for conservation laws, the calculus of variations, gradient flows, among others.
For example, for a system of conservation laws, the existence of a strictly convex entropy
for the system ensures its hyperbolicity; the relative entropy methods via convexity
have proved very useful
for establishing the existence, stability, and structure of solutions, as well as hydrodynamic
limits of large particle systems, for various PDEs.
The notions of convexity appropriate for multi-dimensional problems,
such as polyconvexity, quasiconvexity, and rank-one convexity, have been fundamental
in the theory of nonlinear PDEs
and the calculus of variations
({\it cf.} \cite{Ball77,Ball80,Dacorogna,Evans1,Evans2,Morrey,Dafermos1,LieroMielke}).

\section{Entropy and Convexity for Nonlinear PDEs and Related Areas}

The topics of this theme issue are cross-disciplinary.
The following areas are brought together in the issue:
hyperbolic conservation laws, elliptitc/parabolic equations,
the calculus of variations, continuum mechanics, kinetic theory,
statistical physics, probability, plasma physics,
astrophysics,
materials science,
dynamical systems, optimal transportation, differential geometry, among others.
The underlying connection between them is entropy and convexity.

\smallskip
The paper by Adams, Dirr, Peletier \& Zimmer \cite{ADPZ} provides a survey of recent activities
in deriving and explaining macroscopic evolution equations via entropy
for multi-particle systems.
Heat flows are gradient flows in the space of probability measures, whose
microscopic origin can be justified by using the formalism
of statistical mechanics.
This is done by computing entropy as a rate function via the theory of large
deviations and by identifying the variational structures concerning fluctuations
around the gradient flow, since the zero entropy path is the most likely path (mean-field)
which is the gradient flow; the other possible paths permitted in the microscopic model
corresponds to fluctuations around such a path with a cost (entropy).
The variational structure can be related to the theory of optimal mass transportation
by using an asymptotic expansion ($\Gamma$-convergence) in a kind of inviscid limit
setting.
The study of such phenomena is advocated at a general level by pulling together
the results from different communities ({\it cf.} PDEs, probability,
statistical mechanics, and gradient flows) and exposing related ideas as a whole.
In particular, the emerging potential of a bigger picture is present,
namely that entropies and thermodynamic free energies can be derived via large deviation
principles, at least for some non-equilibrium situations, and
the two concepts of large deviation principles for stochastic particle systems and gradient
flows are closely entwined.
The approach advocated in  \cite{ADPZ} is different from the established hydrodynamic
limit passage and extends a link that is well known in the equilibrium situation.

Blesgen \& Chenchiah in \cite{BlesgenChenchiah} consider the effects of elastic
energy ({\it i.e.} stress) on the evolution of microstructure formation on smaller length scales,
known as coarsening, which are modelled by time-dependent nonlinear PDEs.
Coarsening is driven by both the interfacial and elastic energy.
The equilibrium state for coarsening driven purely by the interfacial energy
is a single spherical inclusion that minimizes the interfacial area.
In contrast, the equilibrium state for coarsening driven by
the elastic energy alone is mathematically challenging and is a microstructure
on an infinitesimal scale.
In \cite{BlesgenChenchiah}, the key experimental observations pertaining to coarsening in elastic solids
are summarised;
the Cahn-Larch\'{e} model, a generalisation of the Cahn-Hillard equation that
incorporates a
quasiconvex
elastic energy density obtained from relaxation, is analyzed;
some recent developments concerning the two-scale model
for elasticity-influenced coarsening are presented,
including several motivations, analytical results,
and comparisons of the model with experimental results; and some further
mathematical reflections and questions are addressed.

The paper by Brenier \cite{Brenier} provides deep, somewhat unexpected, connections of convexity and entropy with the mathematical
theory of convection through the mathematical concept of rearrangement.
Rearrangement theory is about reorganising a given function or map in some specific order,
for example, as an increasing function
or a map with convex potential.
The convection of a fluid leads to fluid parcels
being continuously reorganised in a stabler way under
the action of the buoyancy force.
Convection is one of the most important phenomena in nature with many important applications.
The connections are exposed through deep insights into relations between the so-called transport-collapse method to solve
conservation laws, some interesting social science models, and the structure of some convection models in meteorological sciences.
In particular, the theory enlightens the mechanism involved in the hydrostatic
limit of the Navier-Stokes Boussinesq equations. The limit is obtained
via a relative entropy method under a natural convexity condition.

Hyperbolic systems of conservation laws are one of the most important classes of nonlinear PDEs.
As a consequence of the Second Laws of Thermodynamics, systems of conservation laws arising
in continuum physics are endowed
with an entropy function of the state variable.
In many cases, the entropy function is convex, which is the case for the Euler equations
for elastic fluids, the prototype of hyperbolic systems of conservation laws.
Convexity of the entropy leads to
the local well-posedness of the Cauchy problem in a Sobolev space of
sufficiently high order, as well as the $L^2$-stability of the solution even within the broader
class of entropy solutions ({\it cf.} \cite{BS-book,Chen,Dafermos-book}).
However,
convexity of the entropy function is not ubiquitous; such examples include
important systems of conservation laws in continuum mechanics, thermomechanics, and electrodynamics,
especially the equations of elastodynamics,
for which convexity of the entropy is incompatible with geometric invariance dictated by physics.
Quite often, failure of convexity of the entropy function is encountered in systems endowed with
involutions.
As observed by Dafermos \cite{Dafermos1}, involutions may compensate for the loss of convexity;
indeed, under the assumption that the contingent entropy function is convex merely in the direction of
the involution cone in state space, it is shown in \cite{Dafermos1} that the Cauchy problem is still
locally well-posed in the class
of classical solutions, and that classical solutions are unique and stable even within the broader
class of entropy solutions.
In the process, technical subtleties that were glossed over in earlier treatments are also presented in detail.
The theory is an important generalization of the classical local well-posedness
theory to nonlinear systems of conservation laws equipped with involutions and partially convex entropies, which
include many important physical systems.

The paper by Evans \cite{Evans3} concerns monotonicity methods for nonlinear PDEs arising from variational problems.
Monotonicity methods and entropy methods are strongly related.
For monotonicity formulas, various expressions involving the solution are integrated over
a ball $B(x, r)$ of centre  $x\in \R^n$ and radius $r$ in the Euclidean space $\R^n$,
in order to get useful differential inequalities
determining how these integrals depend on the radius $r$.
The artistry for this approach (as well as the entropy approach) is  in the design of the
precise expressions that are integrated.
To illustrate the approach, it is shown in \cite{Evans3} how certain explicit integral formulas
are derived involving the solutions of elliptic PDEs corresponding to stationary points
of functionals having the form
$$
\int_U F(Du)dx \qquad \mbox{for}\,\,  u : U\subset \R^n\to \R^m,
$$
where $Du$ is the gradient of $u$.
Three
important special cases are treated as illustrations:
(i) $m = 1$ and $F(p) = |p|$;
(ii) $m > 1$ and $F(M) =\frac{1}{2} |M|^2$;
and (iii) $m = n$ and $F(M) =\frac{1}{2}|M|^2 +\frac{1}{(\det M)^\gamma}$.
The approach and ideas presented will stimulate further work
in this direction.

Liero \& Mielke in \cite{LieroMielke} analyze systems of parabolic PDEs involving diffusion, drift, and reaction.
As observed in  \cite{LieroMielke}, under natural physical assumptions,
these systems have the structure of a gradient flow with respect to
an entropy functional and dissipation Riemannian metric given in terms of a so-called Onsager operator, which
is a sum of a diffusion part of Wasserstein type and a reaction part.
New methods are provided for establishing geodesic $\lambda$--convexity of the entropy
functional by purely differential methods, circumventing arguments from mass transportation.
Contraction properties of these gradient flows with respect to the intrinsic dissipation
metric are analyzed.
Sufficient conditions for contractivity are derived for several important physical examples
including a drift-diffusion system, which provides a survey on the applicability of the theory.

The paper by Liu \cite{Liu} samples some of recent analytical developments
in the study of dissipation
around the entropy methods
for hyperbolic conservation laws, viscous conservation laws, and the Boltzmann equation.
As discussed in \cite{Liu} through different examples,
several types of dissipation, such as the viscosity and heat conductivity,
the nonlinearity, and the coupling of distinct
characteristics, occur in both the hyperbolic and viscous systems of conservation laws.
In particular,
the relationship between kinetic
theory and compressible fluid dynamics is explained,
and dissipation caused by the intermolecular collisions in kinetic theory
is addressed.
In addition, the importance of dissipation
due to the boundary effects is emphasized.
The Green's function approach, a qualitative way based on concrete constructions,
is another useful
approach to understand the dissipation and the relation
between gas dynamics and kinetic theory.
There is much room for future progress for
the entropy methods, the Green's function approach,
and other possibilities in this direction.

Penrose in \cite{Penrose} provides a clear illustration
of the importance of phase space volume in the definition
of entropy and effective irreversible behaviour for chaotic dynamical systems.
More precisely, a new non-equilibrium entropy or trajectory entropy
for chaotic dynamical systems is proposed, whose main feature is that it
does not require any use of the concept of a macroscopic state of the dynamical
system. For any given $\varepsilon>0$ and two phase points lying on the two endpoints
of a given trajectory of the system, the so-called dynamical self-correlation of
the two endpoints is roughly the conditional probability that, if the system is
started at initial time from a phase point chosen uniformly in a ball, then its
phase point at the terminal time will lie in an $\varepsilon$-neighbourhood
of the original terminal phase point.
The proposed entropy is then the logarithm of the inverse
of the dynamical self-correlation. The main issue is the conjecture
made
that the dynamical self-correlation converges in the limit of diverging
time (between the chosen phase points of the trajectory) to the quotient of the
phase space measure of the $\varepsilon$-ball and the total phase space measure
of the given system.
Partial support for this conjecture is provided through
two examples of dynamical systems.
The conjecture allows one
to connect the limit of the self-correlation to the equilibrium entropy, and thus
the definition of the new entropy becomes a natural object to consider, that is,
the non-equilibrium entropy is proportional to the logarithm of the quotient
of the phase space measure of the $\varepsilon$-ball and the self-correlation,
for the phase points.
The idea is illustrated by using an example based on Arnold's `cat' map, which also
demonstrates that it is possible to have irreversible behaviour, involving a large increase of entropy,
in a chaotic system with only two degrees of freedom.

Saint-Raymond in \cite{Saint-Raymond}
discusses Boltzmann's $H$--theorem and its essential role in rigorous derivation
of fluid dynamics equations from the Boltzmann equation, which justify
the mathematical significance of the notion of entropy.
The main focus is on the study of hydrodynamic limits in the framework of solutions defined by the physical
energy and entropy bounds, at both the kinetic and fluid level.
The analogies between several kinetic and fluid models are detailed by examining
the use of the entropy inequality and the implications for functional analysis.
The modulated entropy method and the kinetic formulation approach are suggested.
The hydrodynamic regime which has been best understood is the one that leads
to the incompressible Navier-Stokes equations, which is indeed the only asymptotics
of the Boltzmann equation for which an optimal convergence result is known via the scaled relative
entropy inequality ({\it cf.} Golse \& Saint-Raymond \cite{GSR}).

In Slemrod \cite{Slemrod}, three notions of admissibility for weak solutions
are discussed through the isentropic Euler equations of gas dynamics:
the viscosity criterion,
the entropy inequality (the thermodynamically admissible isentropic solutions),
and the viscosity-capillarity criterion.
The Chapman-Enskog expansion for the Boltzmann equation, when truncated to orders beyond
the first order (Burnett and super-Burnett approximations),
is known to produce approximations of
the compressible Euler equations, for which the equilibrium flow is unstable or at most conditionally
stable. The dominant view is that the problem is not the Chapman-Enskog expansion but the
truncation of the expansion.
The focus of this paper is on a simplified model for a one-dimensional flow of Grad's linearized thirteen
moment system. For this system, the Chapman-Enskog expansion can be exactly summed, and
the dispersion relation is calculated accurately to all orders.
The novelty is that the
induced entropy inequality is computed, which has two additional terms: one
corresponds to the capillarity contribution to the energy and the other to the viscous dissipation.
Recent results of DeLellis \& Szekelyhidi \cite{DS} have shown that the usual concept of
entropy solutions is inadequate to imply uniqueness for the incompressible
and the compressible Euler equations. This raises the issue of how this `paradox' is interpreted
or alternatively
how the concept of entropy weak solutions is improved/replaced by an alternative concept.
It is conjectured in \cite{Slemrod} that the energy dissipation, obtained when the compressible Euler
equations are viewed from the perspective of kinetic theory, has additional terms relative to the
one obtained when viewed as a limit from the compressible Navier-Stokes system.
 This is indeed
justified through the simplified model.
Further investigations and explorations are needed to understand the non-uniqueness issue,
which is closely related to entropy and convexity.

\section{Outlooks}

As we have briefly discussed above, entropy and convexity have played
an important role in the recent developments in the analysis of nonlinear PDEs and related areas.
The results in the papers of this theme issue present a deeper understanding
of existing nonlinear methods via entropy and convexity
and their underlying connections in different
areas of PDEs,
and open up challenging new research problems and interdisciplinary research directions.
Further developments of unifying nonlinear methods and ideas via entropy and
convexity will be useful for solving some longstanding, challenging problems
in nonlinear conservation laws, the calculus of variations, and other areas.

In particular, since the seminal work of Jacques Hadamard in the early 1920's concerning
the classification of linear PDEs, the research
community in PDEs has been largely partitioned by the approaches taken to
the mathematical analysis of different classes of PDEs (hyperbolic, parabolic,
and elliptic). However, advances in mathematical research on PDEs over the
last 30 years have been making it increasingly clear that many difficult
questions faced by the community are at the boundaries of this classification
or, indeed, go beyond this narrow classification. Several
important sets of nonlinear PDEs that arise in fluid mechanics, solid mechanics,
materials science, conservation laws, the calculus of variations, geometry, among many others,
involve PDEs of mixed
type ({\it e.g.} mixed hyperbolic-elliptic type, mixed hyperbolic-parabolic type).
The unification of the mathematical theory for disparate classes of nonlinear
PDEs is at the cutting edge of modern mathematical research and has
important implications not only for mathematicians, but also for the wider scientific
communities who use nonlinear PDEs.
Moreover, a particular challenge concerns the understanding of the roles
played in multi-dimensional
dynamics by convexity conditions such as quasiconvexity which are central
to the corresponding variational problems
for statics.
It is our belief that nonlinear methods and ideas via entropy and convexity will also
play a fundamental role in the analysis of nonlinear PDEs of mixed type and for solving
longstanding and newly emerging fundamental open problems
involving such PDEs, which deserve our special attention.

\bigskip
{\bf Acknowledgements}.
We thank the Royal Society for the support of the International Scientific Seminars
on \textit{``Entropy and Convexity for Nonlinear Partial Differential Equations''},
held at the Kavli Royal Society International Centre
for the Advancement of Science at Chicheley Hall, June 16--17, 2010,
which was the origin of this Theme Issue.
J. M. Ball was supported in part by the UK EPSRC Science and Innovation
award to the Oxford Centre for Nonlinear PDE (EP/E035027/1),
by the European Research Council under the European Union's Seventh Framework
Programme (FP7/2007-2013)/ERC grant agreement no. 291053, and
by a Royal Society--Wolfson Research Merit Award.
G.-Q. Chen was supported in part by the UK EPSRC Science and Innovation
award to the Oxford Centre for Nonlinear PDE (EP/E035027/1),
and a Royal Society--Wolfson Research Merit Award.

\end{document}